\documentclass[12pt]{amsart}
\usepackage{amsmath,amssymb, fixltx2e, enumerate} 
\usepackage{mathtools}
\usepackage{units}

\begin{document}
\long\def\/*#1*/{}

\bibliographystyle{plain}
\newtheorem{thm}{Theorem}[section]
\newtheorem{lem}[thm]{Lemma}
\newtheorem{dfn}[thm]{Definition}
\newtheorem{cor}[thm]{Corollary}
\newtheorem{prop}[thm]{Proposition}
\newtheorem{clm}[thm]{Claim}
\theoremstyle{remark}
\newtheorem{exm}[thm]{Example}
\newtheorem{rem}[thm]{Remark}
\def\N{{\mathbb N}}
\def\G{{\mathbb G}}
\def\Q{{\mathbb Q}}
\def\R{{\mathbb R}}
\def\C{{\mathbb C}}
\def\P{{\mathbb P}}
\def\Z{{\mathbb Z}}
\def\v{{\mathbf v}}
\def\x{{\mathbf x}}
\def\O{{\mathcal O}}
\def\M{{\mathcal M}}

\def\E{{\mathbb E}}
\def\kbar{{\bar{k}}}
\def\tr{\mbox{Tr}}
\def\id{\mbox{id}}
\def\qed{{\tiny $\clubsuit$ \normalsize}}

\newcommand{\Qbar}{\overline{k}}
\newcommand{\bpf}{\noindent{\it Proof:}}
\renewcommand{\theenumi}{\alph{enumi}}

\title{Distribution of cokernels of ($n$+$u$) $\times$ $n$ matrices over $\Z_p$}

\author{Ling-Sang Tse}

\begin{abstract}
Let $n, u \geq 0$, $M$ be a ($n$+$u$) $\times$ $n$ matrices over $\Z_p$, and $G$ be a finite abelian p-group group. We find that the probability that the cokernel of $M$ is isomorphic to $\Z_p^u \oplus G$ as $n$ goes to infinity is exactly what is expected from Cohen-Lenstra heuristics for the classical case when $u$ is negative.
\end{abstract}

\maketitle

\vspace{-5pt}


\section{Outline of the paper}
\noindent
The goal of this paper is to prove the following theorem: 
\begin{thm}\label{thm:mainthm}

Let $n, u \geq 0$, $M$ be a randomly chosen ($n$+$u$) $\times$ $n$ matrix over $\Z_p$, and $G$ be a finite abelian p-group. Then
$$ \lim_{n \rightarrow \infty} \P \, (\textnormal {coker } M \cong \Z_p^u \oplus G) = \frac{\prod_{i=1}^{\infty} (1-p^{-i-u})}{|G|^u |\textnormal {Aut } G|}. $$

\end{thm}

\noindent
{\it Notation}: For the rest of this paper, let $n, u \geq 0$, $M$ be a randomly chosen ($n$+$u$) $\times$ $n$ matrix over $\Z_p$, and $G$ be a finite abelian p-group group. 
Also, let $\mu$ be the standard Haar measure on the set of ($n$+$u$) $\times$ $n$ matrices over $\Z_p$ such that 
$\mu ($\{($n$+$u$) $\times$ $n$ matrices over $\Z_p$\}) = 1. \\

\noindent
We prove Theorem \ref{thm:mainthm} using the following lemmas: 

\begin{lem}\label{lem: rank} 
For any $(n + u) \times n$ matrix $M$ over $\Z_p$, the probability that $M$ has rank $n$ is 1. 
\end{lem}

\begin{lem} \label{lem: moment1}
If M is a random $(n + u) \times n$ over $\Z_p$, then 

$$\lim_{n \rightarrow \infty} \, \E_{\mu}\,  (\textnormal{\#Sur } (\textnormal{coker }M, G)) = |G|^u$$
\end{lem}

\noindent
Note that from Lemma \ref{lem: rank}, the probability that $M$ has cokernel isomorphic to $\Z_p^u \oplus T$ for some abelian p-group $T$ is 1, 
so what the lemma really meant to say is $\lim_{n \rightarrow \infty} \, \E_{\mu}\,  (\#\textnormal{Sur } (\Z_p^u \oplus T, G)) = |G|^u$.

\begin{lem}\label{lem: moment2}
If M is a random $(n + u) \times n$ over $\Z_p$, then 

$$\lim_{n \rightarrow \infty} \, \E_{\mu}\,  (\textnormal{\#Sur } ( (\textnormal{coker }M) [p^\infty], G)) = |G|^u$$

\end{lem}

\noindent
Similarly as in Lemma \ref{lem: moment1}, what Lemma \ref{lem: moment2} meant is $\lim_{n \rightarrow \infty} \, \E_{\mu}\,  (\textnormal{\#Sur } (T, G)) = |G|^u$.

%
%
\section{Proofs}
\noindent
{\it Proof of Lemma \ref{lem: rank}}. 
An alternate, more combinatorial, proof of this lemma is presented in the appendix in {\it A.1}. Here, we present a much nicer proof.\\

An $(n+u) \times n$ matrix has rank n if and only if it it has an $n \times n$ submatrix with rank n, so it is suffucient to show that the probability that an $n \times n$  matrix has rank $n$ is 1.

Claim \ref{lem: rank}.1: Fix a polynomial $f$ in $\Z_p$ over n variables. The probability of picking a random point $(a_1, ..., a_n)$ in $(\Z_p)^n$ so that $f$ vanishes at $(a_1, ..., a_n)$ is 0.

Proof of Claim \ref{lem: rank}.1 by induction:
If n is 1, this is clear because the polynomial has finitely many roots.

Assume this is true for $n$ - 1 variables, and let $f$ be a polynomial in n variables over $\Z_p$. At any given given point in $(a_1, ..., a_n) \in (Z_p)^n$, consider $g(x_n) := f(a_1, ..., a_{n-1}, x_n)$ as a polynomial in the one variable $x_n$. The coefficients of $g(x_n)$ are polynomials in $n$ - 1 variables, so by induction, the probability of picking a random point $(a_1, ..., a_n)$ in $(\Z_p)^n$ so that all of the coefficients of the corresponding polynomial $g(x_n)$ are 0 is 0. Thus, the probability that $g(x_n)$ is a non-trivial polynomial is 1, so as in the 1 variable case, the probability that $g(a_n) = f(a_1, .., a_n)$ = 0 is 0. This completes the proof of Claim \ref{lem: rank}.1.

An $n \times n$ matrix has rank n if and only if its determinant is 0, and the determinant of an $n \times n$ matrix is a fixed polynomial in $n^2$ variables evaluated at a random point in $(\Z_p)^{n^2}$. Thus, by Claim \ref{lem: rank}.1, the probability that an $n \times n$ matrix has rank $n$ is 0. \\

%
\noindent
{\it Proof of Lemma \ref{lem: moment1}}. 
Let $M$ be a random $(n + u) \times n$ matrix over $\Z_p$. We have
\begin{align*}
     \E_{\mu}\,  (\textnormal{Sur } (\textnormal{coker }M, G)) &= \int\limits_{\{ M: \Z_p^n \rightarrow \Z_p^{n+u} \}}  \sum\limits_{\{ \phi: \textnormal{ coker } M \twoheadrightarrow \, G \}} 1 \\
     &= \sum\limits_{\{ \phi: \, \Z_p^{n+u} \twoheadrightarrow \, G \}} 
     \int\limits_{\begin{Bsmallmatrix}
 M: \, \Z_p^n \rightarrow \Z_p^{n+u} 
\\ \textnormal{such that im }  M \subseteq \textnormal{ ker} \phi
\end{Bsmallmatrix} }
  d\mu \\
 &= \sum\limits_{\{ \phi: Z_p^{n+u} \twoheadrightarrow G \}} |G|^{-n}
  \end{align*}
  
The last equality follows because ker $\phi$ has index $|G|$ in $\Z_p^{n+u}$, so the probability that any of the basis element maps into ker $\phi$ is $|G|^{-1}$. $M$ is determined by where the n basis elements are mapped to, so the probability that im $M \subseteq $ ker $\phi$ is $|G|^{-n}$. 

Therefore, we have
\begin{align*}
   \lim_{n \rightarrow \infty} \, \E_{\mu}\,  \#(\textnormal{Sur } (\textnormal{coker }M, G))  &= \lim_{n \rightarrow \infty} \, \sum\limits_{\{ \phi: Z_p^{n+u} \twoheadrightarrow G \}} |G|^{-n} \\
   &= \lim_{n \rightarrow \infty} \, \sum\limits_{\{ \phi: Z_p^{n+u} \rightarrow G \}} |G|^{-n} \\
   &= \frac{|G|^{n+u}}{|G|^n}\\
   &= |G|^u   
\end{align*}

The second last equality follows because as $n \rightarrow \infty$, $\#$Hom ($\Z_p^{n+u}, G) \sim \#$Sur ($\Z_p^{n+u}, G)$ (We present a proof of this in the appendix in {\it A.2}). 

This completes the proof of lemma \ref{lem: moment1}.\\

\noindent
{\it Proof of Lemma \ref{lem: moment2}}. 

As in the note beneath the statement of Lemma \ref{lem: moment1}, let $\Z_p^u \oplus T$ be the cokernel of an $(n+u) \times n$ matrix. We have

\begin{align*}
    \sum\limits_{H \leq G} \# \textnormal{Sur } (\Z_p^u \oplus T, H) &= \# \textnormal{Hom } (\Z_p^u \oplus T, G) \\
    &= \#\textnormal{Hom } (\Z_p^u , G) \#\textnormal{Hom } (T , G) \\
    &= |G|^u \big{(} \sum\limits_{H \leq G} \# \textnormal{Sur } (T, H) \Big{)} 
\end{align*}

Therefore, 

$$\# \textnormal{Sur } (T, G) = \sum\limits_{H \leq G} \frac{\# \textnormal{Sur } (\Z_p^u \oplus T, H)}{|G|^u} - \sum\limits_{H \lneq G} \# \textnormal{Sur } (T, H)$$

We prove Lemma \ref{lem: moment2} by inducting on the order of $|G|$. 

If $G$ is the trivial group, then it is clear that $\E_{\mu}\,  (\textnormal{\#Sur } ( (\textnormal{coker }M) [p^\infty]) = |G|^u =1 $.

We proceed with the induction, and assume that the hypothesis is true for any group of order less than $|G|$.  

\begin{align*}
    \E_{\mu}\,  (\textnormal{\#Sur } (T, G)) &= \sum\limits_{H \leq G} \frac{ \E_{\mu} (\# \textnormal{Sur } (\Z_p^u \oplus T, H))}{|G|^u} - \sum\limits_{H \lneq G} \E_{\mu} (\# \textnormal{Sur } (T, H)) \\
    &= \sum\limits_{H \leq G} \frac{|H|^u}{|G|^u} - \sum\limits_{H \lneq G} |H|^{-u}
\end{align*}

where the second equality follows by Lemma \ref{lem: moment1} and by the induction hypothesis.

\noindent
Claim \ref{lem: moment2}.1: If $n \, \big{|} \, |G|$, 

$$\# \begin{Bmatrix}
 \textnormal{subgroups of } G
\\ \textnormal{of order }  n
\end{Bmatrix}  =
\# \begin{Bmatrix}
 \textnormal{subgroups of } G
\\ \textnormal{of index }  n
\end{Bmatrix} $$

\noindent
Proof of Claim \ref{lem: moment2}.1: We show

$$
    \#\begin{Bmatrix}
 \textnormal{subgroups of } G
\\ \textnormal{of order }  n
\end{Bmatrix} = \# \begin{Bmatrix}
 \textnormal{subgroups of } \hat{G}
\\ \textnormal{of index }  n
\end{Bmatrix}   =
\#\begin{Bmatrix}
 \textnormal{subgroups of } G
\\ \textnormal{of index }  n
\end{Bmatrix}
$$


To show the second equality, we may define an isomorphism $\phi_2$ from $G$ to $\hat{G}$: For a cyclic group $G_0$ with generator $\alpha_0$, we map $\alpha_0$ to the character $\big{(} \alpha_0 \mapsto \textnormal{e}^{\frac{2\pi i}{|G|}})$. Then for a general finite abelian group $G$, we decompose $G$ into a product of cyclic groups and extend the map the natural way by taking the product of the maps on the cyclic groups. It is easy to check that $\phi_2$ is a well-defined isomorphism, so the second equality follows.

To show the first equality, we define a bijective map of sets $\phi_1$ mapping a subgroup of $G$ of order n to a subgroup of $\hat{G}$ of index $n$. Let $H^\perp := \{ \chi \in \hat{G} \,\, | \,\, \chi (h) = 1 \textnormal{ for all } h \in H \}$, and we define $\phi_1$ by $\phi_1(H) = H^\perp$. By the universal property of quotients, given a character $\chi \in \hat{G}$, $H$ is contained n the kernel of $\chi$ if and only if there exists a homomorphism $\bar{\chi}: \nicefrac{G}{H} \rightarrow \C$ such that $\chi (g) = \bar{\chi} (\bar{g})$ for all $g \in G$. Thus,  every $\chi \in H^\perp$ corresponds to a character $\bar{\chi} \in \widehat{\nicefrac{G}{H}}$ by reduction mod $H$, and so $H^\perp \, \cong \,$  $\widehat{\nicefrac{G}{H}}$. As in the previous paragraph, $\widehat{\nicefrac{G}{H}} \, \cong \, \nicefrac{G}{H}$, and it is easy to check that $H^\perp$ is a subgroup of $\hat{G}$, so $H^\perp$ is a subgroup of $\hat{G}$ of index n. Thus $\phi_1$ is well-defined and is injective.

It remains to check that $\phi_2$ is a bijection: By definition of $(H^\perp)^\perp$,  $H \subseteq (H^\perp)^\perp$ (where we treat $h \in H \subseteq \hat{\hat{G}}$ as evaluation at $h$). Since $|(H^\perp)^\perp| = |H|$, we have $(H^\perp)^\perp = H$. 

The proof of Claim \ref{lem: moment2}.1 is complete. \\

Therefore, 
\begin{align*}
    \E_{\mu}\,  (\textnormal{\#Sur } (T, G)) &= \sum\limits_{H \leq G} \frac{|H|^u}{|G|^u} - \sum\limits_{H \lneq G} |H|^{-u} \\
    &= |G|^{-u}
\end{align*}

This completes the proof of Lemma \ref{lem: moment2}, and we can finally prove our main theorem.\\

%
%
\noindent
{\it Proof of Theorem \ref{thm:mainthm}}. 
We apply the following theorem, due to Wood:

\begin{thm}\cite[Theorem 8.3]{Wood} \label{thm: wood}
Let $X_n$ be a sequence of random variables taking values in finitely generated abelian groups. Let $a$ be a positive integer and $A$ be the set of (isomorphism classes of)
abelian groups with exponent dividing a. Suppose that for every $G \in A$, we have
$$\lim\limits_{n\rightarrow \infty}
\E(\# \textnormal{Sur }(X_n, G)) \leq | \wedge^2 G|.$$
Then for every H ∈ A, the limit $\lim\limits_{n \rightarrow\infty} \P(X_n \otimes \nicefrac{\Z}{a\Z} \cong H)$ exists, and for all $G \in A$ we
have
$$
\sum\limits_{H \in A}
\lim\limits_{n \rightarrow \infty}
\P(X_n \otimes \Z/a\Z \cong H) \# \textnormal{ Sur}(H, G) \leq |\wedge^2 G |$$
Suppose $Y_n$ is a sequence of random variables taking values in finitely generated abelian
groups such that for every $G \in A$, we have
$$\lim\limits_{n \rightarrow \infty}
\E(\# \textnormal{ Sur}(Y_n, G)) \leq | \wedge^2 G |.$$
Then, we have that for every $H \in A$
$$\lim_{n\rightarrow \infty}
\P(X_n \otimes \Z/a\Z \cong H) = \lim\limits_{n\rightarrow \infty}
\P(Y_n \otimes \Z/a\Z \cong H).$$
\end{thm}

Pick random finite p-groups $Y_n$ with probability $$\frac{\prod\limits_{i=1}^\infty (1-p^{-i-u}) } {|G|^u |\textnormal{ Aut } G|},$$ and denote the expectation values with respect to this probability measure $\E_Y$. 

By Theorem 3.20 in [Majumder], 
$$ \sum\limits_{G \textnormal{ is a p-group}} \frac{1} {|G|^u |\textnormal{ Aut } G|} = \prod\limits_{i=1}^\infty (1-p^{-i-u})^{-1}, $$
so this is a well-defined probability distribution. 

We apply the following theorem, also due to Wood, to get
$$\lim\limits_{n \rightarrow \infty} \E_{\mu} (\# \textnormal{Sur } (T, G)) = \lim\limits_{n \rightarrow \infty} \E_Y( \# \textnormal{Sur} (Y, G)) = |G|^{-u}$$

\begin{thm} \cite[Theorem 3.2]{Wood2} \label{thm: wood2}
Let $Y$ be a random p-group chosen with probability $$\frac{\prod\limits_{i=1}^\infty (1-p^{-i-u}) } {|G|^u |\textnormal{ Aut } G|}.$$ Then for every finite abelian group $G$ with exponent dividing $a$, we have
$$\E(\# \textnormal{Sur}(Y, G)) = |G|^{-u}
.$$
\end{thm}

Let $a = p|G|$, and we apply Theorem \ref{thm: wood} to get
$$  \lim_{n \rightarrow \infty} \P_{\mu} \, ( T \cong G) = \frac{\prod_{i=1}^{\infty} (1-p^{-i-u})}{|G|^{u} |\textnormal {Aut } G|}.
$$ and so
\begin{align*}
   \lim_{n \rightarrow \infty} \P_{\mu} \, ( \textnormal {coker } M \cong \Z_p^u \oplus G) &= \frac{\prod_{i=1}^{\infty} (1-p^{-i-u})}{|G|^{u} |\textnormal {Aut } G|}.
\end{align*}
This completes the proof of Theorem \ref{thm:mainthm}. Similarly, we may get the following corollary for a matrix over $\Z$:
\begin{cor} \label{cor: corollary}
Let $n, u \geq 0$, $M$ be a randomly chosen ($n$+$u$) $\times$ $n$ matrix over $\Z$, $G$ be a finite abelian group, and $P$ be the set of primes dividing $|G|$. Then
$$ \lim_{n \rightarrow \infty} \P \, ( \textnormal {coker } M \cong \Z^u \oplus G) = \prod\limits_{p \in P} \frac{\prod_{i=1}^{\infty} (1-p^{-i-u})}{|G|^u |\textnormal {Aut } G|}. $$
\end{cor}

\section{Appendix}
\noindent
{\it A.1 Alternative proof of Lemma \ref{lem: rank}}. 

\noindent
We prove Lemma \ref{lem: rank} by showing the following:

1. The rank of $M$ is equal to n if and only if $\exists$ $e > 0$ such that ($p^e \Z)^{n} \subset$ im $M$. To be precise, 
what is meant by ($p^e \Z)^{n} \subset$ im $M$ in the rest of this proof is whether there exists an embedding of ($p^e \Z)^{n}$ into im $M$. 

Then by Nakayama's lemma, for all $e \geq 0, e' > e$, $\P$ (rank $M$ = n) $\geq \P$ \big{(}$(\frac{p^e\Z}{p^{e'}\Z})^{n} \subset $ im $M$ (mod $p^{e'}$)$\big{)}$. 

2. If $e' > e > 0, \, \P \big{(}(\frac{p^e\Z}{p^{e'}\Z})^{n} \subset $ im $M$ $\big{)}$ $\ge$ $(\frac{(n+u)!}{u!})^{1/n}p^{-ue-u-e-1}$, which goes to 0 as $e$ goes to infinity. \\

\noindent
\underline{Step 1:} We show that the rank of $M$ is equal to n $\iff \exists$ $e > 0$ such that ($p^e \Z)^{n} \subset$ im $M$.\\

$(\Longrightarrow)$
\begin{align*}
& \mbox{rank } M = n \\
& \iff \exists \mbox { an } n \times n \mbox { matrix of minor with nonvanishing determinant} \\
& \iff \mbox { im } M \mbox { is minimallly generated by } \{v_i\}_{i = 1}^{n} \mbox { for some } v_i = (v_{ij})_{1 \leq j \leq n} \in \Z_p^n
 \end{align*}
 

Let $e = \underset{1 \leq i,j \leq n}{\mbox{max}} \{$ ord$_P$ ($v_{ij}$) $\}$. Then 

$$ \textnormal {im } M \mbox{ is minimallly generated by } \{ v_i \}_{i=1}^{n} \iff (p^e\Z_p)^{n} \subset \mbox { im } M $$ 

$( \Longleftarrow ) $
Let $e_i$ be the standard i$\textsuperscript {th}$ basis vector. We may assume that there exists $e > 0$ such that $p^ee_i \in $ im $M$ for all $i$ satisfying 1 $\leq i \leq n$. Then for all $i$, there exist $v_1, ..., v_n \in \Z_p^n$ 
such that $Mv_i = p^ee_i$.  Then $M(p^{-e}v_i) = e_i$ for all $i$ (if we extend the action of $M$ to $\Q_p^n$), so there exists an $n \times n$ matrix 
of minors of $M$ such that det $M \neq 0$. It follows that the rank $M$ is $n$. \qed \\

\underline{Step 2:} We finish the proof by showing that if $e' > e > 0, \, $ then $$ \P \big{(}(\frac{p^e\Z}{p^{e'}\Z})^{n} \subset \textnormal{ im } M \big{)} \ge p^{-e(n-1)n + (n+u)^2}\prod\limits_{i=1}^{n+u}(1-p^{-i}).$$

Let $\bar{M} \in \Big{(} \frac{\Z}{p^{e'}\Z} \big{)}^{(n+u) \times n}$, and suppose $\Big{(} \frac{p^e\Z}{p^{e'}\Z} \big{)}^n \not\subseteq$ im $\bar{M}$. Then we may post-compose $M$ with an automorphism of $\Big{(} \frac{p^e\Z}{p^{e'}\Z} \big{)}^{(n+u) \times n}$ fixing $\Big{(} \frac{p^{e+1}\Z}{p^{e'}\Z} \big{)}^{(n+u) \times n}$ such that im $\bar{M}$ does not contain any of the vectors in the set $\{ (0, ..., 0, p^e, 0, ..., 0) \}_{i = n-1}^{n+u}$, with the $p^e$ in the $i^{th}$ place. 

Counting the number of choices for $\phi \circ \bar{M}$ is equivalent to counting the number of $(n+u) \times n$ matrices of the form 

$$
\begin{bmatrix}
 \frac{\Z}{p^{e'}\Z} ... ... ...\frac{\Z}{p^{e'}\Z}
\\ \\ \\  \frac{\Z}{p^{e'}\Z} ... ... ... \frac{\Z}{p^{e'}\Z} \\
\frac{\Z}{p^{e'-e-1}\Z} ... ... ...\frac{\Z}{p^{e'-e-1}\Z} \\ \\ \\ 
\frac{\Z}{p^{e'-e-1}\Z} ... ... ...\frac{\Z}{p^{e'-e-1}\Z}
\end{bmatrix} ,
$$

with $n-1$ rows of $\frac{\Z}{p^{e'}\Z}$ at the top and $u+1$ rows of $\frac{\Z}{p^{e'-e-1}\Z}$. There are $p^{n((n-1)e' + (n+1)(e'-e-1))}$ of them. Choosing $\phi$ is equivalent to choosing an automorphism of $\Big{(} \frac{\Z}{p\Z} \Big{)}^{n+u},$ and there are 
$p^{(n+u)^2}\prod\limits_{i=1}^{n+u} (1-p^{-i})$ of them. 

There are $p^{n(n+u)e'}$ hoomomorphisms from $\Big{(} \frac{\Z}{p^{e'}\Z} \Big{)}^{n}$ to $\Big{(} \frac{\Z}{p^{e'}\Z} \Big{)}^{n+u}$ in total. 

Therefore, 

\begin{align*}
      \P \big{(}(\frac{p^e\Z}{p^{e'}\Z})^{n} \not\subset \textnormal{ im } M \big{)}&\leq 
      \frac{p^{n((n-1)e'+(u+1)(e'-e-1))}p^{(n+u)^2} \prod\limits_{i=1}^{n+u}(1-p^{-1})}{p^{n(n+u)e'}} 
     \\ &= p^{-e(n-1)n + (n+u)^2}\prod\limits_{i=1}^{n+u}(1-p^{-i}), 
\end{align*}

which goes to 0 as e goes to $\infty$. \\

Therefore, for all $e \geq 0, e' > e$, 

\begin{align*}
    \P ( \textnormal{ rank } M = n) \geq \P \Big{(}  \big{(} \frac{p^e\Z}{p^{e'}\Z} \big{)}^n\subseteq \textnormal{ im } M  \Big{)},  
\end{align*}

which goes to 1 as e goes to $\infty$. \\

\noindent
{\it A.2 Proof of $\#$Hom ($\Z_p^{n+u}, G) \sim \#$Sur ($\Z_p^{n+u}, G)$:} 

We have 
\begin{align*}
    \# \textnormal{ Sur } (\Z_p^{n+u}, G) &= \# \{ \bar{\phi} : \big{(} \frac{\Z}{p\Z} \big{)}^{n+u} \twoheadrightarrow \frac{G}{pG} \} \cdot \{ \phi: \Z_p^{n+u} \rightarrow G \mid \bar{\phi} = 0\} \\
    &= (p^{n+u} - 1)(p^{n+u}-p)... (p^{n+u}-p^{k-1}) \# \{ \phi: \Z_p^{n+u} \rightarrow G \mid \bar{\phi} = 0\} \\
    &= p^{(n+u)k}(1-p^{-(n+u)})... (1-p^{-(n+u)+k-1})\# \{ \phi: \Z_p^{n+u} \rightarrow G \mid \bar{\phi} = 0\} \\
    &= \# \textnormal{ Hom } (\Z_p^{n+u}, G)(1-p^{-(n+u)})(1-p^{-(n+u)+1})... (1-p^{-(n+u)+k-1}),
\end{align*}

which goes to $\#$Hom ($\Z_p^{n+u}, G$) as $n$ goes to $\infty$.

\section*{Acknowledgements}
I would like to thank my PhD advisor Jacob Tsimerman for his insightful comments and helpful discussions in writing this paper.

\end{document}